\documentclass[letterpaper, 10 pt, conference]{ieeeconf}  

\usepackage{cite}
\usepackage{amsmath,amssymb,amsfonts}
\usepackage{tikz,tikz-cd}

\IEEEoverridecommandlockouts                              

\title{\LARGE \bf
Geometric Integrators for Nonholonomic Systems on Lie Groups
}

\author{Viyom Vivek, David Mart{\'i}n de Diego, and Ravi N. Banavar, \IEEEmembership{Senior Member, IEEE}
\thanks{}
\thanks{V. Vivek is with the Centre for Systems \& Control, Indian Institute of Technology Bombay, India (e-mail: viyomvivek@iitb.ac.in).}
\thanks{D. M. de Diego is with Instituto de Ciencias Matem{\'a}ticas Madrid, Spain (e-mail: david.martin@icmat.es).}
\thanks{R. N. Banavar is with 
the Centre for Systems \& Control, Indian Institute of Technology Bombay, India (e-mail: banavar@iitb.ac.in).}
}

\begin{document}

\maketitle
\thispagestyle{empty}
\pagestyle{empty}

\begin{abstract}

We present a general framework for constructing structure-preserving numerical integrators for nonholonomically constrained mechanical systems evolving on Lie groups using retraction maps. Retraction maps generalize the exponential map and provide a convenient tool for performing numerical integration on manifolds. In nonholonomic mechanics, the constraints restrict the dynamics to a nonintegrable distribution rather than the entire tangent bundle. Using the Hamel formulation, the equations of motion can be expressed in local coordinates adapted to this constraint distribution. We then specialize the framework to the case of Lie groups, where both the dynamics and the constraints exhibit symmetries, allowing a simplified formulation of the numerical scheme. The resulting integrator respects the constraint distribution and enforces the nonholonomic constraints at each discrete time step. The approach is illustrated using the Suslov problem.

\end{abstract}

\section{Introduction}
\label{sec: introduction}
This paper develops a natural extension of the retraction-map–based framework introduced in \cite{Vivek2025lcss} for constructing structure-preserving numerical integrators for mechanical systems evolving on Lie groups. In the present work, we adapt this framework to incorporate nonholonomic constraints, extending the work of \cite{Vivek2025ecc}, and developing the first of its kind intrinsic integrator for a nonholonomically constrained mechanical system on a Lie group.

Classical examples of nonholonomically constrained mechanical systems include the rolling motion of a wheel without slipping and the motion of a sleigh subject to a knife-edge constraint; see, for example, \cite{Bloch2015}. Such constraints restrict the admissible velocities of the system at every point in the configuration space, but cannot, in general, be integrated to constraints on the configuration variables alone. Consequently, they define nonintegrable distributions on the configuration manifold rather than lower-dimensional configuration submanifolds; see \cite{Neimark1972} for further discussion.

Mechanical systems subject to nonholonomic constraints exhibit a remarkably rich geometric structure involving constrained variational principles, distributions on Lie groups, and reduced dynamics on Lie algebras. These geometric features play a central role in the development of numerical schemes that faithfully capture the qualitative aspects of the continuous dynamics. In particular, preserving the constraint distribution and the Lie group structure in the discretization is crucial for obtaining integrators with favorable long-term behavior. In this work, we exploit these geometric properties to construct numerical integrators that remain compatible with both the constraint and the Lie group geometry of the configuration space.

For background on mechanics on Lie groups we refer to \cite{Abraham2008}, while standard references on nonholonomic mechanics include \cite{Bloch2015}. An overview of retraction maps and their role in the construction of geometric numerical integrators can be found in \cite{Absil2009, Barbero-Linan2023}. We also refer to \cite{Celledoni2019, Cortes-Martinez2001, Ferraro2008, Mclachlan2006} for representative numerical approaches to nonholonomic systems.

\textbf{Outline of the paper.} We begin with a physically intuitive introduction to nonholonomically constrained mechanical systems, followed by a precise geometric formulation of their dynamics and its specialization to the setting of Lie groups. We then review the retraction-map–based framework for numerical integration on Lie groups. Building on this framework, we develop a modification that incorporates the nonholonomic constraints and ensures that they are satisfied at each discrete time step. Finally, we illustrate the resulting method with a representative example that demonstrates the applicability of the proposed approach.
\section{Nonholonomic Mechanics}
\label{sec: nonholonomic mechanics}
\subsection{Physical Motivation}
\label{subsec: physical motivation}
\begin{center}
\begin{tikzpicture}[scale = 0.8]
        \draw[gray,fill=gray,fill opacity=0.1] (1,3) to[out=15,in=105] (3,1) to[out=-5,in=-175] (7,1) to[out=105,in=15] (5,3) to[out=-175,in=-5] (1,3);
        \node at (5,0.5) {\scriptsize $Q = \mathbb{R}^2 \times \mathbb{S}^1 \times \mathbb{S}^1$};
        \draw[fill] (4,2.3) circle (1pt);
        \node at (4.9,2.05) {\tiny $q = (x,y,\theta,\phi)$};
        \draw[dotted] (4,2.3)--(4,3.8);
        \draw[cyan,fill=cyan,fill opacity=0.1] (3.5,3.3)--(5.5,3.3)--(4.5,4.3)--(2.5,4.3)--cycle;
        \draw[thick,blue] (3.1,3.7)--(4.9,3.9);
        \node at (2.8,4.6) {\scriptsize $\textcolor{cyan}{T_qQ}$};
        \node at (2,4.65) {\scriptsize $\textcolor{cyan}{\mathbb{R}^4 \cong}$};
        \node at (2.7,3.5) {\scriptsize $\textcolor{blue}{D_q}$};
        \node at (2.05,3.55) {\scriptsize $\textcolor{blue}{\mathbb{R}^2 \cong}$};
        \draw[ultra thick,blue,-stealth] (4,3.8)--(4.6,3.87);
        \node at (5.1,4.15) {\tiny $\textcolor{blue}{\dot{q}=(\dot{x},\dot{y},\dot{\theta},\dot{\phi})}$};
        \draw[dotted] (-0.5,1.5)--(-2,1);
        \draw[-stealth,blue,thick] (-0.5,1.5)--(0.2,1.5);
        \node at (0.3,1.6) {\tiny $\textcolor{blue}{\dot{x}}$};
        \draw[-stealth,blue,thick] (-0.5,1.5)--(0,1.833);
        \node at (0.1,1.9) {\tiny $\textcolor{blue}{\dot{y}}$};
        \draw[-stealth,blue,thick] (-0.5,1.5)--(-1,1.333);
        \node at (-1.1,1.6) {\tiny $\textcolor{blue}{r\dot{\phi}}$};
        \draw[white,fill=gray,fill opacity=0.2] (-3,1)--(0.5,1)--(2,2)--(-1.5,2)--cycle;
        \draw[thick,fill=cyan,fill opacity=1] (-0.5,3) ellipse (0.5 and 1.5);
        \draw[fill] (-0.5,3) circle (1pt);
        \draw[->] (-3,1)--(1,1);
        \draw[->] (-3,1)--(-1.2,2.2);
        \draw[->] (-1,1) arc (0:35:0.5);
        \draw[fill] (-0.5,1.5) circle (1pt);
        \draw[dotted] (-0.5,1.5)--(-0.5,3);
        \draw[dotted] (-0.5,3)--(0,2.2);
        \draw[<-] (-0.5,2.4) arc (-90:-50:0.5);
        \node at (-0.8,1.15) {\tiny $\theta$};
        \node at (-0.3,2.2) {\tiny $\phi$};
        \node at (-0.2,1.3) {\tiny $(x,y)$};
        \node at (-0.7,2.5) {\tiny $r$};
\end{tikzpicture}
\end{center}
To motivate the notion of nonholonomic constraints, consider the motion of a wheel, or equivalently a vertical disk of radius $r$, rolling on a planar surface without slipping. The configuration space of this system is $Q = \mathbb{R}^2 \times \mathbb{S}^1 \times \mathbb{S}^1$ with configuration variables $q := (x,y,\theta,\phi)$ as illustrated above. Here $(x,y)$ denotes the position of the point of contact of the wheel with the plane, $\theta$ represents the orientation of the wheel in the plane, and $\phi$ denotes the angle of rotation of the wheel about its axis. 

The velocity of the system at configuration $q$ is $\dot{q} := (\dot{x},\dot{y},\dot{\theta},\dot{\phi}) \in T_q Q \cong \mathbb{R}^4$. The rolling-without-slipping constraint requires that the point of contact of the wheel with the plane remain instantaneously at rest. This leads to the velocity constraints
$\dot{x} = r \dot{\phi}\cos\theta$ and 
$\dot{y} = r \dot{\phi}\sin\theta$. Equivalently, these constraints can be written in matrix form as
$A(q)\dot{q} = 0$,
where
$A(q) := \left(
\begin{smallmatrix}
1 & 0 & 0 & -r\cos\theta \\
0 & 1 & 0 & -r\sin\theta
\end{smallmatrix}
\right)$. Thus, at each configuration $q \in Q$, the admissible velocities are restricted to a two-dimensional subspace $D_q \cong \mathbb{R}^2$ of $T_qQ \cong \mathbb{R}^4$. In other words, the constraints define a distribution $D \subset TQ$ on the configuration manifold $Q$ that specifies the allowable directions of motion.

A key feature of these constraints is that they cannot be integrated to constraints depending solely on the configuration variables. In particular, the condition $A(q)\dot{q} = 0$ does not imply the existence of a function $B : Q \to \mathbb{R}$ such that the constraints can be written in the form $dB(q)\dot{q} = 0$. 
Constraints of this type, which restrict the admissible velocities without corresponding to constraints on the configuration variables alone, are called nonholonomic constraints. 

It is important to note that the forces acting at the point of contact, which enforce the constraint, perform no work on the system. Equivalently, the corresponding constraint forces are orthogonal to the admissible velocities in $D$. This property is fundamental in the formulation of the dynamics via the Lagrange-d'Alembert principle.
\subsection{Geometric Description}
\label{subsec: geometric description}
We now give a geometric description of nonholonomic constraints. Let $Q$ be the configuration manifold of a mechanical system with $n$ degrees of freedom, with local coordinates $(q^i)$ for $i \in \{1,\dots,n\}$. The space of configurations and velocities is the tangent bundle $TQ$, equipped with the induced coordinates $(q^i,\dot{q}^i)$. Similarly, the space of configurations and momenta is the cotangent bundle $T^*Q$, with the corresponding dual coordinates $(q^i,p_i)$. 

The presence of nonholonomic constraints restricts the admissible velocities to a subbundle (or distribution) $D \subset TQ$. For each configuration $q \in Q$, the fiber $D_q \subset T_qQ$ specifies the set of all admissible velocities at that configuration. Locally, a distribution of rank $m \leq n$ can be specified by system of linear velocity constraints of the form $A^\alpha_i(q) \dot{q}^i = 0$, for $\alpha \in \{m+1,\dots,n\}$.

In the nonholonomic setting, this distribution is generally nonintegrable and therefore does not arise as the tangent bundle of a submanifold of $Q$. Consequently, $D$ represents the space of configurations together with their admissible velocities.

The annihilator of this distribution, denoted $D^\circ \subset T^*Q$, is a codistribution consisting of covectors that vanish on $D$. Locally, $D^\circ$ is spanned by the constraint one-forms $A^\alpha := A^\alpha_i(q) dq^i$. From a physical perspective, the fiber $D_q^\circ \subset T_q^*Q$ represents the space of constraint forces enforcing the nonholonomic constraints at the configuration $q \in Q$.

The equations of motion governed by the Lagrange-d'Alembert principle are
\begin{align} 
    \frac{d}{dt} \left( \frac{\partial \mathcal{L}}{\partial \dot{q}^i} \right) - \frac{\partial \mathcal{L}}{\partial q^i} = \lambda_\alpha A_i^\alpha \quad \text{and} \quad
    A_i^\alpha \dot{q}^i = 0 \label{Euler-Lagrange nonholonomic}
\end{align}
where $\mathcal{L}: TQ \to \mathbb{R}$ is the Lagrangian of the system and $\lambda_\alpha$ are Lagrange multipliers enforcing the nonholonomic constraints, see \cite{Bloch2015}. The first set of equations describes the dynamics in the presence of constraint forces, while the second set imposes a set of kinematic constraints.

Using the Legendre transform $\mathbb{F} \mathcal{L} : TQ \to T^*Q$, we define the Hamiltonian $\mathcal{H} : T^*Q \to \mathbb{R}$ of the system by
\begin{align*}
    \mathcal{H}(q^i,p_i) := \langle p_i, \dot{q}^i \rangle - \mathcal{L}(q^i,\dot{q}^i) \quad \text{where} \quad p_i := \frac{\partial \mathcal{L}}{\partial \dot{q}^i}.
\end{align*}
The equations of motion \eqref{Euler-Lagrange nonholonomic} in Hamiltonian form are
\begin{align}
    \dot{q}^i = \frac{\partial \mathcal{H}}{\partial p_i}, \quad \dot{p}_i = -\frac{\partial \mathcal{H}}{\partial q^i} + \lambda_\alpha A_i^\alpha \quad \text{and} \quad A_i^\alpha \frac{\partial \mathcal{H}}{\partial p_i} = 0 \label{Hamilton nonholonomic}
\end{align}
where the last equation expresses the nonholonomic constraint in Hamiltonian variables, see \cite{Bloch2015}. Here the multipliers $\lambda_\alpha$ represent the constraint forces enforcing the admissibility of the velocities. The presence of these forces causes the dynamics to deviate from the standard Hamiltonian flow and, in particular, leads to a violation of symplecticity, see \cite{Grabowski2009}.
\subsection{Hamel's Formulation}
\label{subsec: hamels formulation}
A choice of local coordinates $(q^i)$ on $Q$ induces a local frame $\left\{ \frac{\partial}{\partial q^i} \right\}$ on $TQ$, with respect to which a velocity vector can be written as $\dot{q}^i \frac{\partial}{\partial q^i}$, where the coefficients $\dot{q}^i$ are the components of the velocity. In the presence of nonholonomic constraints, it is often convenient to work with a local frame $\{X_i(q)\}$ adapted to the constraint distribution. Using the Riemannian metric $\mathcal{G}$ on $Q$, typically induced by the kinetic energy of the system, the tangent bundle can be decomposed as $TQ = D \oplus D^{\perp,\mathcal{G}}$. This decomposition yields an adapted frame $\{X_a(q), X_\alpha(q)\}$, where $X_a(q)$ span the constraint distribution $D$ and $X_\alpha(q)$ span its $\mathcal{G}$-orthogonal complement $D^{\perp,\mathcal{G}}$, with indices $a \in \{1,\dots,m\}$. A velocity vector can then be expressed in the form $y^a X_a(q) + y^\alpha X_\alpha(q)$, where the coefficients $(y^a,y^\alpha)$ are called quasivelocities, \cite{Bloch2009}. This representation of the dynamics in terms of an adapted frame and quasivelocities is known as the Hamel formulation.

Define the adapted Lagrangian $\mathcal{L} : D \oplus D^{\perp,\mathcal{G}} \to \mathbb{R}$ which in local coordinates takes the form $\mathcal{L}(q^i,y^a,y^\alpha)$. We then define the restricted Lagrangian $l : D \to \mathbb{R}$ by restricting $\mathcal{L}$ to the constraint distribution, that is, $l := \left. \mathcal{L} \right|_D$, which locally depends only on $(q^i,y^a)$. 

In terms of these variables, the nonholonomic equations of motion \eqref{Euler-Lagrange nonholonomic} take the form
\begin{align}
    \frac{d}{dt} \left( \frac{\partial l}{\partial y^a} \right) + C_{ab}^c y^b \frac{\partial l}{\partial y^c} - X^i_a \frac{\partial l}{\partial q^i} = 0 \,\, \text{and} \,\, \dot{q}^i = X^i_a y^a. \label{Euler-Lagrange Hamel}
\end{align}
Here the coefficients $C_{ab}^c$ arise from the Lie brackets of the adapted frame and are defined by $[X_a,X_b] = C_{ab}^c X_c + C_{ab}^\alpha X_\alpha$, where $b,c \in \{1,\dots,m\}$. These coefficients are known as the structure functions of the adapted frame, see \cite{Bloch2009} for further details.

Using the Legendre transform $\mathbb{F}l : D \to D^*$, we define the restricted Hamiltonian $h : D^* \to \mathbb{R}$ by
\begin{align*}
    h(q^i,\rho_a) := \langle \rho_a, y^a \rangle - l(q^i,y^a) \quad \text{where} \quad \rho_a := \frac{\partial l}{\partial y^a}.
\end{align*}
The equations of motion \eqref{Euler-Lagrange Hamel} in Hamiltonian form are
\begin{align}
    \dot{q}^i = X^i_b \frac{\partial h}{\partial \rho_b} \quad \text{and} \quad \dot{\rho}_a = -C_{ab}^c \rho_c\frac{\partial h}{\partial \rho_b} - X^i_a \frac{\partial h}{\partial q^i}. \label{Hamilton Hamel}
\end{align}
A more intrinsic formulation of these equations, expressed in terms of an almost Poisson structure on $D^*$, can be found in \cite{Balseiro2009, Yoshimura2007}. We now specialize the preceding construction to Lie groups.
\subsection{Lie Group Setting}
\label{subsec: on lie groups}
We now consider mechanical systems with nonholonomic constraints evolving on a Lie group $G$. A left-invariant nonholonomic constraint distribution is completely determined by a subspace $\mathfrak{d} \subset \mathfrak{g}$ of the Lie algebra $\mathfrak{g}$, which is typically not a Lie subalgebra. The entire distribution can be reconstructed by left translation, namely $D_g := (L_g)_* \mathfrak{d}$ for every $g \in G$, where $L_g$ denotes left multiplication on the group. Similarly, a left-invariant mechanical system is completely described by a reduced Lagrangian $\mathcal{L} : \mathfrak{g} \to \mathbb{R}$, with the dynamics governed by Euler-Poincar{\'e} equations; see \cite{Marsden2013}. 

Using the Hamel formulation, we obtain the decomposition $\mathfrak{g} = \mathfrak{d} \oplus \mathfrak{d}^\perp$, where the orthogonal complement is taken with respect to the inner product that induces the left-invariant Riemannian metric representing the kinetic energy of the system. Choosing a basis $\{e_a,e_\alpha\}$ of $\mathfrak{g}$ adapted to this decomposition, that is, $\mathfrak{d} = span \{e_a\}$ and $\mathfrak{d}^\perp = span \{e_\alpha\}$, any element $\xi \in \mathfrak{g}$ can be written as $\xi = \xi^a e_a + \xi^\alpha e_\alpha$, where the coefficients $(\xi^a,\xi^\alpha)$ are the quasivelocities. 

Define the restricted reduced Lagrangian $l: \mathfrak{d} \to \mathbb{R}$ by restricting $\mathcal{L}$ to the constraint subspace, namel $l := \left. \mathcal{L} \right|_{\mathfrak{d}}$, which locally depends only on the variable $\xi^a$. The reduced nonholonomic equations of motion in Hamel form are then given by
\begin{align}
    \frac{d}{dt} \left( \frac{\partial l}{\partial \xi^a} \right) + C_{ab}^c \xi^b \frac{\partial l}{\partial \xi^c} = 0 \quad \text{and} \quad \dot{g} = g (\xi^a e_a). \label{Euler-Poincare Hamel}
\end{align}
where the second equation is the reconstruction equation used to recover the dynamics on the group $G$. Here the coefficients $C_{ab}^c$ are the structure constants of the Lie algebra $\mathfrak{g}$ with respect to the adapted basis, defined by $[e_a,e_b] = C_{ab}^c e_c + C_{ab}^\alpha e_\alpha$, and they capture the noncommutativity of the Lie group $G$.

Taking the Legendre transform $\mathbb{F}l : \mathfrak{d} \to \mathfrak{d}^*$, we define the restricted reduced Hamiltonian $h : \mathfrak{d}^* \to \mathbb{R}$ by 
\begin{align*}
    h(\mu_a) := \langle \mu_a, \xi^a \rangle - l(\xi^a) \quad \text{where} \quad \mu_a := \frac{\partial l}{\partial \xi^a}.
\end{align*}
The equations of motion \eqref{Euler-Poincare Hamel} in Hamiltonian form are
\begin{align}
    \dot{g} = g \left( \frac{\partial h}{\partial \mu_b} e_b \right) \quad \text{and} \quad \dot{\mu}_a = -C_{ab}^c \mu_c \frac{\partial h}{\partial \mu_b}. \label{Lie-Poisson Hamel}
\end{align}
The first equation is the reconstruction equation on the Lie group, while the second describes the reduced dynamics on $\mathfrak{d}^*$. A more intrinsic formulation of these equations, expressed in terms of almost Lie-Poisson structure, can be found in \cite{Balseiro2009, Yoshimura2007}.
\section{Retraction-Map-Based Framework}
\label{sec: retraction map based framework}
\subsection{Retraction Maps}
\label{subsec: retraction maps}
A smooth map $\mathcal{R} :  G \times \mathfrak{g} \to G$ satisfying (i) $\mathcal{R}(g,0) = g$ and (ii) $\left. \frac{d}{dt} \right|_{t=0} \mathcal{R}(g,t\xi) = T_eL_g(\xi)$ for every $g \in G$ and $\xi \in \mathfrak{g}$ is called a (left-trivialized) retraction map on the Lie group $G$. 
\begin{center}
    \begin{tikzpicture}[scale = 0.6]
        \draw[gray,fill=gray,fill opacity=0.1] (1,3) to[out=15,in=105] (3,1) to[out=-5,in=-175] (9,1) to[out=105,in=15] (7,3) to[out=-175,in=-5] (1,3);
        \node at (2,1.5) {\scriptsize $G$};
        \draw[fill] (3.5,2.5) circle (1pt);
        \draw (3.5,2.5) to[out=-25,in=115] (4.5,1.5);
        \draw[fill] (4.2,2) circle (1pt);
        \node at (3.2,2.5) {\tiny $e$};
        \node at (4.5,1.2) {\tiny $\tau(t\xi)$};
        \node at (4.7,2.1) {\tiny $\tau(\xi)$};
        \draw[fill] (7,2) circle (1pt);
        \draw (7,2) to[out=-25,in=115] (8,1);
        \draw[fill] (7.7,1.5) circle (1pt);
        \node at (6.7,2) {\tiny $g$};
        \node at (9,0.6) {\scriptsize $\mathcal{R}(g,t\xi):=L_g(\tau(t\xi))$};
        \node at (10,1.5) {\scriptsize $\mathcal{R}(g,\xi):=L_g(\tau(\xi))$};
        \draw[dotted] (3.5,2.5)--(3.5,4);
        \draw[fill=cyan,fill opacity=0.1,cyan] (3,3.5)--(5.5,3.5)--(4.5,4.5)--(2,4.5)--cycle;
        \draw[cyan,fill] (3.5,4) circle (1pt);
        \draw[cyan,thick,-stealth] (3.5,4)--(4.3,3.7); 
        \node at (3.2,4) {\tiny $\textcolor{cyan}{0}$};
        \node at (4,4.1) {\tiny $\textcolor{cyan}{\xi}$};
        \node at (2,3.8) {\scriptsize $\textcolor{cyan}{\mathfrak{g}}$};
        \draw[dotted] (7,2)--(7,3.7);
        \draw[gray,fill=gray,fill opacity=0.3] (6.5,3.2)--(9,3.2)--(8,4.2)--(5.5,4.2)--cycle;
        \draw[fill] (7,3.7) circle (1pt);
        \draw[thick,-stealth] (7,3.7)--(7.8,3.4);
        \node at (6.7,3.7) {\tiny $0$};
        \node at (7.8,3.8) {\tiny $T_eL_g(\xi)$};
        \node at (9,4) {\scriptsize $T_gG$};
\end{tikzpicture}
\end{center}
Such a map can be constructed as follows. Let $\tau : \mathfrak{g} \to G$ be a local diffeomorphism satisfying (i) $\tau(0) = e$ and (ii) $\left. \frac{d}{dt} \right|_{t=0} \tau(t\xi) = \xi$
for all $\xi \in \mathfrak{g}$. Then the map defined by $\mathcal{R}(g,\xi) := L_g(\tau(\xi))=g\tau(\xi)$ is a retraction map as shown above. The most common example of such a map $\tau$ is the Lie group exponential map. However, other useful choices exist for particular Lie groups, depending on computational convenience or numerical properties; see \cite{Bou-Rabee2009} for further examples and discussion.

We now introduce the notion of logarithmic derivatives, which will be needed in the construction of the integrator, see \cite{Michor2008}. Let $\tau : \mathfrak{g} \to G$ be a local diffeomorphism as above. Consider the affine line $\xi + t\eta$ in $\mathfrak{g}$. Under $\tau$, this line gets mapped to a curve in $G$. The velocity of this curve at $t=0$ can be translated back to the identity by left translation, as illustrated below.
\begin{center}
\begin{tikzpicture}[scale = 0.8]
        \draw[gray,fill=gray,fill opacity=0.1] (1,3) to[out=15,in=105] (3,1) to[out=-5,in=-175] (9,1) to[out=105,in=15] (7,3) to[out=-175,in=-5] (1,3);
        \node at (2,1.5) {\scriptsize $G$};
        \draw[->] (4.7,3.2) to[out=-75,in=75] (4.7,2.3);
        \node at (5,2.5) {\scriptsize $\tau$};
        \draw[fill] (4.5,2) circle (1pt);
        \draw (3,2.1) to[out=5,in=135] (6,1.1);
        \node at (3,0.5) {\tiny $R_{\tau(\xi)^{-1}}(\tau(\xi+t\eta))$};
        \draw (3.1,1.1) to[out=45,in=180] (5.9,2.2);
        \node at (6,0.5) {\tiny $L_{\tau(\xi)^{-1}}(\tau(\xi+t\eta))$};
        \node at (4.5,1.7) {\tiny $e$};
        \draw[fill] (7.55,2) circle (1pt);
        \node at (7.2,1.7) {\tiny $\tau(\xi)$};
        \draw (6.5,2.7) to[out=-25,in=115] (8.1,1.2);
        \node at (8.4,0.7) {\tiny$\tau(\xi+t\eta)$};
        \draw[dotted] (4.5,2)--(4.5,4);
        \draw[fill=cyan,fill opacity=0.1,cyan] (3.5,3.3)--(7,3.3)--(6,4.7)--(2.5,4.7)--cycle;
        \draw[cyan,thick,-stealth] (4.5,4)--(5.25,4); 
        \node at (5.3,4.2) {\tiny \textcolor{cyan}{$\xi$}};
        \draw[cyan,thick,-stealth] (5.25,4)--(5.7,3.82);
        \node at (5.8,4) {\tiny \textcolor{cyan}{$\eta$}};
        \draw[cyan] (4,4.5)--(6.5,3.5);
        \node at (7,3.15) {\tiny $\textcolor{cyan}{\xi+t\eta}$};
        \node at (4.3,4.1) {\tiny $\textcolor{cyan}{0}$};
        \node at (2.5,3.8) {\scriptsize $\textcolor{cyan}{\mathfrak{g}}$};
        \draw[thick,blue,-stealth] (4.5,4)--(5.4,3.4);
        \node at (5.5,3) {\scriptsize \textcolor{blue}{$d^L_\xi\tau(\eta)$}};
        \draw[thick,blue,-stealth] (4.5,4)--(3.6,3.4);
        \node at (3.5,3) {\scriptsize \textcolor{blue}{$d^R_\xi\tau(\eta)$}};
        \draw[cyan,fill] (4.5,4) circle (1pt);
\end{tikzpicture}
\end{center}
This motivates the definition of a smooth map $d^L_\xi\tau : \mathfrak{g} \to \mathfrak{g}$, given by $d^L_\xi \tau(\eta) := \left. \frac{d}{dt} \right|_{t=0} \tau(\xi)^{-1}\tau(\xi+t\eta)$ for all $\eta \in \mathfrak{g}$. This map is called the left logarithmic derivative of $\tau$ at $\xi \in \mathfrak{g}$. The right one is defined analogously. The two are related by the adjoint action as $d^R_\xi \tau = Ad_{\tau(\xi)} \circ d^L_\xi \tau$. 
\subsection{Discretization Maps}
\label{subsec: discretization maps}
A smooth map $\mathcal{D} : G \times \mathfrak{g} \to G \times G$, which assigns $(g,\xi) \mapsto (\mathcal{D}_1(g,\xi), \mathcal{D}_2(g,\xi))$, is called a (left-trivialized) discretization map on the Lie group $G$ if it satisfies (i) $\mathcal{D}_1(g,0) =\mathcal{D}_2(g,0)=g$ and (ii) $\left. \frac{d}{dt} \right|_{t=0} \mathcal{D}_2(g,t\xi) - \left. \frac{d}{dt} \right|_{t=0} \mathcal{D}_1(g,t\xi) = T_eL_g(\xi)$ for every $g \in G$ and $\xi \in \mathfrak{g}$. Such a map has two important properties as illustrated below.

First, a retraction $\mathcal{R} : G \times \mathfrak{g} \to G$ gives rise to a family of discretization maps $\mathcal{D} : G \times \mathfrak{g} \to G \times G$ defined as $\mathcal{D}(g,\xi) := (\mathcal{R}(g,-s\xi),\mathcal{R}(g,(1-s)\xi))$ for each $s \in [0,1]$.

Second, a discretization $\mathcal{D}: G \times \mathfrak{g} \to G \times G$ is locally invertible in a neighborhood of $(g,0) \in G \times \mathfrak{g}$. This follows directly from the inverse function theorem.
\begin{center}
    \begin{tikzpicture}[scale = 0.8]
        \draw[gray,fill=gray,fill opacity=0.1] (1,3) to[out=15,in=105] (3,1) to[out=-5,in=-175] (9,1) to[out=105,in=15] (7,3) to[out=-175,in=-5] (1,3);
        \node at (2,1.5) {\scriptsize $G$};
        \draw[fill] (3.5,2.5) circle (1pt);
        \draw (3.5,2.5) to[out=-25,in=115] (4.5,1.5);
        \draw[fill] (4.2,2) circle (1pt);
        \node at (3.2,2.5) {\tiny $e$};
        \node at (4.5,1.2) {\tiny $\tau(t\xi)$};
        \node at (4.7,2.1) {\tiny $\tau(\xi)$};
        \draw[fill] (7,2) circle (1pt);
        \draw (7,2) to[out=-25,in=115] (8,1);
        \draw (7,2) to[out=155,in=-5] (5.5,2.3);
        \draw[fill] (7.7,1.5) circle (1pt);
        \draw[fill] (6,2.25) circle (1pt);
        \node at (6.8,1.8) {\tiny $g$};
        \node at (8,0.6) {\tiny $L_g(\tau(t\xi))$};
        \node at (9.5,1.9) {\scriptsize $\mathcal{D}_2(g,\xi):=L_g(\tau((1-s)\xi))$};
        \node at (7.5,2.5) {\scriptsize $\mathcal{D}_1(g,\xi):=L_g(\tau(-s\xi))$};
        \draw[dotted] (3.5,2.5)--(3.5,4);
        \draw[fill=cyan,fill opacity=0.1,cyan] (3,3.5)--(5.5,3.5)--(4.5,4.5)--(2,4.5)--cycle;
        \draw[cyan,fill] (3.5,4) circle (1pt);
        \draw[cyan,thick,-stealth] (3.5,4)--(4.3,3.7); 
        \node at (3.2,4) {\tiny $\textcolor{cyan}{0}$};
        \node at (4,4.1) {\tiny $\textcolor{cyan}{\xi}$};
        \node at (2,3.8) {\tiny $\textcolor{cyan}{\mathfrak{g}}$};
        \draw[dotted] (7,2)--(7,3.7);
        \draw[gray,fill=gray,fill opacity=0.3] (6.5,3.2)--(9,3.2)--(8,4.2)--(5.5,4.2)--cycle;
        \draw[fill] (7,3.7) circle (1pt);
        \draw[thick,-stealth] (7,3.7)--(7.8,3.4);
        \node at (6.7,3.7) {\tiny $0$};
        \node at (7.8,3.8) {\tiny $T_eL_g(\xi)$};
        \node at (9,4) {\scriptsize $T_gG$};
\end{tikzpicture}
\end{center}
\subsection{Numerical Integrators}
\label{subsec: numerical integrators}
We now highlight the core idea of this approach to constructing numerical integrators on Lie groups. Let $X : G \to TG$ be a vector field whose flow we wish to discretize. 
\begin{center}
    \begin{tikzcd}
        & G \times \mathfrak{g} & & G \times G \arrow[r,"\mathcal{D}^{-1}"] \arrow[d,"\mathcal{D}^{-1}",swap] & G \times \mathfrak{g} \\
        G \arrow[r,"X",swap] \arrow[ur,"\widetilde{X}"] & TG \arrow[u,"TL_{g^{-1}}",swap] & & G \times \mathfrak{g} \arrow[r,"\mathfrak{pr}_1",swap] & G \arrow[u,"\widetilde{X}",swap] 
    \end{tikzcd}
\end{center}
\noindent \textbf{Step 1:} Choose a retraction $\mathcal{R} : G \times \mathfrak{g} \to G$ and a value of the parameter $s \in [0,1]$. Construct the discretization $\mathcal{D} : G \times \mathfrak{g} \to G \times G$ by $\mathcal{D}(g,\xi) := (\mathcal{R}(g,-s\xi),\mathcal{R}(g,(1-s)\xi))$.

\noindent \textbf{Step 2:} Left trivialize the vector field as illustrated in the diagram above (left). Define the map $\widetilde{X} : G \to G \times \mathfrak{g}$ by $\widetilde{X}(g) := TL_{g^{-1}}(X(g))$ for all $g \in G$.

\noindent \textbf{Step 3:} Construct the numerical integrator
\begin{align}
    t \widetilde{X}(\mathfrak{pr}_{1}(\mathcal{D}^{-1}(g_{(k)},g_{(k+1)}))) = \mathcal{D}^{-1}(g_{(k)},g_{(k+1)}) \label{integrator on G}
\end{align}
using the commutative diagram shown above (right).

Here $t \in \mathbb{R}$ is the time step. Its role is to scale the vector field $\widetilde{X}$ so that its values lie within the domain where the discretization map $\mathcal{D}$ is locally invertible, ensuring that the above equation defining the integrator is well posed.
\subsection{Trivialization on Lie Groups}
\label{subsec: trivialization on lie groups}
For analyzing systems that exhibit symmetries in the form of left-invariance, it is particularly useful to work with trivializations. The tangent bundle $TG$ can be trivialized as $TG \cong G \times \mathfrak{g}$ using left translation as alluded to earlier. The key is to perform this trivialization in a way that preserves the Lie group structure of the tangent bundle, see \cite{Esen2014, Grabowska2016}. Similarly, the cotangent bundle $T^*G$ can be trivialized as $T^*G \cong G \times \mathfrak{g}^*$ in an analogous manner. In this case, the goal is not only to preserve the natural Lie group structure of the cotangent bundle but also the canonical symplectic structure; again see \cite{Esen2014, Grabowska2016}.

As a mechanical system evolves either on the tangent bundle (Lagrangian viewpoint) or on the cotangent bundle (Hamiltonian viewpoint), we will trivialize the corresponding second-order bundles in a similar manner. After performing these trivializations, we obtain the following identifications: 
\begin{align*}
    &TTG \cong G \times \mathfrak{g} \times \mathfrak{g} \times \mathfrak{g}, \quad TT^*G \cong G \times \mathfrak{g}^* \times \mathfrak{g} \times \mathfrak{g}^*, \\
    &T^*TG \cong G \times \mathfrak{g} \times \mathfrak{g}^* \times \mathfrak{g}^*, \quad T^*T^*G \cong G \times \mathfrak{g}^* \times \mathfrak{g}^* \times \mathfrak{g}.
\end{align*}
The trivialized version of all maps will henceforth be denoted by placing a tilde over the corresponding symbol, as was done earlier for vector fields.
\subsection{Lifted Numerical Integrators}
\label{subsec: lifted numerical integrators}
For constructing numerical integrators for mechanical systems, we need to lift discretization maps from the Lie group to its tangent or cotangent bundle, see \cite{Barbero-Linan2023}. For Lie groups, these lifted discretization maps can then be trivialized as discussed above; the details can be found in \cite{Vivek2025lcss}.
\subsubsection{Tangent-Lifted}
\label{subsubsec: tangent lifted}
The (left-trivialized) tangent-lifted discretization map can be defined using the following commutative diagram
\begin{center}
    \begin{tikzcd}
        G \times \mathfrak{g} \times \mathfrak{g} \times \mathfrak{g} \arrow[r,"\widetilde{\mathcal{D}}^T"] \arrow[d,shift left,"\widetilde{\kappa}_G"] & G \times \mathfrak{g} \times G \times \mathfrak{g} \arrow[d,shift left,"\mathfrak{sw}_{23}"] \\
        G \times \mathfrak{g} \times \mathfrak{g} \times \mathfrak{g} \arrow[d,"\mathfrak{pr}_{12}",swap] \arrow[r,"\widetilde{T\mathcal{D}}"] \arrow[u,shift left,"\widetilde{\kappa}_G"] & G \times G \times \mathfrak{g} \times \mathfrak{g} \arrow[d,"\mathfrak{pr}_{12}"] \arrow[u,shift left,"\mathfrak{sw}_{23}"] \\
        G \times \mathfrak{g} \arrow[r,"\mathcal{D}",swap] & G \times G
    \end{tikzcd}
\end{center}
\begin{align}
    \widetilde{\mathcal{D}}^T := sw_{23} \circ \widetilde{T\mathcal{D}} \circ \widetilde{\kappa}_G \label{tangent lifted discretization}
\end{align}
where $\kappa_G : TTG \to TTG$ is the canonical involution, see \cite{Michor2008}. The map $\mathfrak{pr}_{12}$ denotes the projection onto the first and second factors, while $\mathfrak{sw}_{23}$ denotes the map that switches the second and third factors.
\subsubsection{Cotangent-lifted}
\label{subsubsec: cotangent lifted}
The (left-trivialized) cotangent-lifted discretization map can be defined using the following commutative diagram
\begin{center}
    \begin{tikzcd}
        G \times \mathfrak{g}^* \times \mathfrak{g} \times \mathfrak{g}^* \arrow[r,"\widetilde{\mathcal{D}}^{T^*}"] \arrow[d,"\widetilde{\alpha}_G",swap] & G \times \mathfrak{g}^* \times G \times \mathfrak{g}^* \arrow[d,"\widetilde{\Phi}"] \\
        G \times \mathfrak{g} \times \mathfrak{g}^* \times \mathfrak{g}^* \arrow[d,"\mathfrak{pr}_{12}",swap]  & G \times G \times \mathfrak{g}^* \times \mathfrak{g}^* \arrow[d,"\mathfrak{pr}_{12}"] \arrow[l,"\widetilde{T^*\mathcal{D}}",swap] \\
        G \times \mathfrak{g} \arrow[r,"\mathcal{D}",swap] & G \times G
    \end{tikzcd}
\end{center}
\begin{align}
    \widetilde{\mathcal{D}}^{T^*} := \widetilde{\Phi}^{-1} \circ (\widetilde{T^*\mathcal{D}})^{-1} \circ \widetilde{\alpha}_G \label{cotangent lifted discretization}
\end{align}
where $\alpha_G : TT^*G \to T^*TG$ is the symplectomorphism dual to the canonical involution, see \cite{Grabowska2016} The map $pr_{12}$ denotes the projection onto the first and second factors. The map $\Phi : T^*G \times T^*G \to T^* (G \times G)$ is the twisted-product symplectomorphism, see \cite{Silva2001}.

We can use tangent-lifted \eqref{tangent lifted discretization} or cotangent-lifted \eqref{cotangent lifted discretization} discretization maps, together with the general recipe presented in \eqref{integrator on G}, to construct second-order numerical integrators. The integrator obtained using \eqref{cotangent lifted discretization} is
\begin{align}
    t \widetilde{X}(\mathfrak{pr}_{12}(&(\widetilde{\mathcal{D}}^{T^*})^{-1}(g_{(k)}, \mu_{(k)},g_{(k+1)},\mu_{(k+1)}))) = \nonumber\\ &(\widetilde{\mathcal{D}}^{T^*})^{-1}(g_{(k)},\mu_{(k)},g_{(k+1)},\mu_{(k+1)}) \label{integrator on T*G}
\end{align}
where the continuous dynamics are governed by the vector field $X: T^*G \to TT^*G$. This integrator preserves symplecticity and will be our primary focus from this point onward.
\section{Nonholonomic Numerical Integrators}
\label{sec: nonholonomic numerical integrator}
Finally, we have assembled all the ingredients needed to construct a numerical integrator for a left-invariant mechanical system on a Lie group subject to left-invariant nonholonomic constraints. The resulting integrator is designed to respect the constraints by preserving the constraint distribution at each discrete time step. 

Let $\iota : D \hookrightarrow TG$ and $\pi : TG \to D$ denote the canonical inclusion and projection maps arising from the Hamel decomposition $TG = D \oplus D^{\perp}$. Recall that, due to left-invariance, $D_g := (L_g)_*\mathfrak{d}$, where $\mathfrak{d} \subset \mathfrak{g}$ and $\mathfrak{g} = \mathfrak{d} \oplus \mathfrak{d}^\perp$. On the dual side, we have the corresponding decomposition $\mathfrak{g}^* = \mathfrak{d}^* \oplus (\mathfrak{d}^\perp)^*$. Using these maps, we constrain the cotangent-lifted discretization map to remain within the constraint codistribution, 
\begin{center}
    \begin{tikzcd}
        G \times \mathfrak{g}^* \times \mathfrak{g} \times \mathfrak{g}^* \arrow[r,"\widetilde{\mathcal{D}}^{T^*}"] & G \times \mathfrak{g}^* \times G \times \mathfrak{g}^* \arrow[d,"\widetilde{\iota^*} \times \widetilde{\iota^*}"] \\
        G \times \mathfrak{d}^* \times \mathfrak{d} \times \mathfrak{d}^* \arrow[r,"\widetilde{\mathcal{D}}^{D^*}",swap] \arrow[u,"\widetilde{T\pi^*}"] & G \times \mathfrak{d}^* \times G \times \mathfrak{d}^*
    \end{tikzcd}
\end{center}
\begin{align}
    \widetilde{\mathcal{D}}^{D^*} := (\widetilde{\iota^*} \times \widetilde{\iota^*}) \circ \widetilde{\mathcal{D}}^{T^*} \circ \widetilde{T\pi^*} \label{discretization on D*}
\end{align}
using the commutative diagram above. Since this constraint-adapted discretization map is constructed as a projected version of a symplectic integrator, it is expected to approximately preserve the energy of the system as well. This is consistent with the physical fact that ideal constraint forces do no work; consequently, the dynamics preserve constant energy surfaces, although they do not preserve phase space volume or symplectic structure, see \cite{Zenkov2003}.

Let $X: G \times \mathfrak{d}^* \to T(G \times \mathfrak{d}^*)$ be the vector field corresponding to the equations of motion \eqref{Lie-Poisson Hamel}. Applying the general recipe presented in \eqref{integrator on G}, we obtain
\begin{align}
     t \widetilde{X}(\mathfrak{pr}_{12}(&(\widetilde{\mathcal{D}}^{D^*})^{-1}(g_{(k)}, \mu_{(k)},g_{(k+1)},\mu_{(k+1)}))) = \nonumber\\ &(\widetilde{\mathcal{D}}^{D^*})^{-1}(g_{(k)},\mu_{(k)},g_{(k+1)},\mu_{(k+1)}) \label{integrator on D*}
\end{align}
which defines a constraint preserving integrator.
\section{Illustrative Example}
\label{sec: illustrative example}
We illustrate the construction of our integrator by applying it to one of the most pedagogically important examples of a left-invariant nonholonomic system evolving on a Lie group, known as the Suslov problem, see \cite{Bloch1996}. 
\subsection{Suslov Problem}
\label{subsec: problem statement}
The Suslov problem describes a rigid body rotating about a fixed point subject to a constraint that removes one component of the body angular velocity. Physically, this corresponds to a mechanism preventing rotation about a chosen body-fixed direction. Its simplicity and geometric clarity make it a standard example for illustrating nonholonomic systems on Lie groups.

Here $G = \mathbb{SO}(3)$, the group of rotation matrices, represents the orientation, while $\mathfrak{g} = \mathfrak{so}(3) \cong \mathbb{R}^3$ represents body angular velocity under the standard identification given by the hat map $\hat{\cdot} : \mathbb{R}^3 \to \mathfrak{so}(3)$, see \cite{Holm2009}. Let $\{e_1,e_2,e_3\}$ be the standard basis of $\mathbb{R}^3$, with respect to which a body angular velocity $\Omega$ can be written as $\Omega^i e_i$, where $i \in \{1,2,3\}$. The kinetic energy is $\frac{1}{2} \mathcal{I}_{ij} \Omega^i \Omega^j$ where $j \in \{1,2,3\}$, and the inertia matrix $\mathcal{I} = \left( \begin{smallmatrix}
    \mathcal{I}_{11} & 0 & 0 \\
    0 & \mathcal{I}_{22} & 0 \\
    0 & 0 & \mathcal{I}_{33}
\end{smallmatrix} \right)$ defines an inner product on $\mathbb{R}^3$. The constraint distribution is obtained by left translation of the subspace $\mathfrak{d} := span \{e_1, e_2\} \cong \mathbb{R}^2$; equivalently, the constraint can be written as $\Omega^3=0$. This constraint is nonholonomic because $\mathfrak{d}$ is not a Lie subalgebra of $\mathfrak{g}$; indeed $[e_1,e_2] = e_3 \notin span\{e_1,e_2\}$. Using the Hamel formulation, we obtain the decomposition $\mathfrak{g} = \mathfrak{d} \oplus \mathfrak{d}^\perp = span \{e_1,e_2\} \oplus span \{e_3\}$. The restricted reduced Lagrangian is $l(\Omega) = \frac{1}{2} \mathcal{I}_{ab} \Omega^a \Omega^b$ where $a,b \in \{1,2\}$. Substituting into the equations of motion \eqref{Euler-Poincare Hamel} yields
\begin{align}
    \mathcal{I}_{11} \dot{\Omega}^1 = 0, \quad \mathcal{I}_{22} \dot{\Omega}^2 = 0, \quad \dot{R} = R \left (\begin{smallmatrix}
        0 & 0 & \Omega^2 \\
        0 & 0 & -\Omega^1 \\
        -\Omega^2 & \Omega^1 & 0
    \end{smallmatrix} \right) \label{Euler Poincare Suslov}
\end{align}
where $R \in \mathbb{SO}(3)$ defines the attitude. Taking the Legendre transform, we obtain the restricted reduced Hamiltonian $h(\Pi) = \frac{1}{2} \mathcal{I}^{ab} \Pi_a \Pi_b$ where $\Pi \in \mathbb{R}^3 \cong \mathfrak{so}(3)^*$ is the body angular momentum and $\mathcal{I}^{ab}$ are the components of the inverse inertia matrix. The equations of motion \eqref{Euler Poincare Suslov} then take the Hamiltonian form
\begin{align}
    \dot{\Pi}_1 = 0, \quad \dot{\Pi}_2 = 0, \quad \dot{R} = R \left( \begin{smallmatrix}
        0 & 0 & \mathcal{I}^{22} \Pi_2 \\
        0 & 0 & -\mathcal{I}^{11} \Pi_1 \\
        -\mathcal{I}^{22} \Pi_2 & \mathcal{I}^{11} \Pi_1 & 0 
    \end{smallmatrix} \right) \label{Lie Poisson Suslov equations}
\end{align}
using \eqref{Lie-Poisson Hamel}. The vector field corresponding to these dynamics is
\begin{align}
    X(R,\Pi_1,\Pi_2) := (&R, \Pi_1,\Pi_2; \nonumber \\
    &R \left( \begin{smallmatrix}
        0 & 0 & \mathcal{I}^{22} \Pi_2 \\
        0 & 0 & -\mathcal{I}^{11} \Pi_1 \\
        -\mathcal{I}^{22} \Pi_2 & \mathcal{I}^{11} \Pi_1 & 0 
    \end{smallmatrix} \right),0,0). \label{Lie Poisson Suslov vector field}
\end{align}
\subsection{Adapted Discretization Map}
\label{subsec: adapted discretization map}
We begin by choosing the retraction $\mathcal{R}(g,\xi) := g \tau(\xi)$ for all $(g,\xi) \in G \times \mathfrak{g}$. Setting the parameter $s=0$ yields the discretization $\mathcal{D}(g,\xi) := (g,g\tau(\xi))$ for all $(g,\xi) \in G \times \mathfrak{g}$. Next, using \eqref{cotangent lifted discretization}, we compute the cotangent-lifted discretization. Its inverse is
\begin{align}
    &(\mathcal{D}^{T^*})^{-1}(g_{(k)},\mu_{(k)},g_{(k+1)},\mu_{(k+1)}) = \nonumber \\
    &(g_{(k)},d^{L*}_{\tau^{-1}(g_{(k)}^{-1}g_{(k+1)})}\tau(\mu_{(k+1)}),\tau^{-1}(g_{(k)}^{-1}g_{(k+1)}), \nonumber \\
    &Ad^*_{(g_{(k)}^{-1}g_{(k+1)})^{-1}}(\mu_{(k+1)})-\mu_{(k)}). \label{cotangent lifted discretization inverse}
\end{align}
Using \eqref{discretization on D*}, we obtain the constraint-adapted discretization map. Its inverse is
\begin{align}
    &(\mathcal{D}^{D^*})^{-1}(g_{(k)},\bar{\mu}_{(k)},g_{(k+1)},\bar{\mu}_{(k+1)}) = \nonumber \\
    &(g_{(k)},d^{L*}_{\bar{\tau}^{-1}(g_{(k)}^{-1}g_{(k+1)})}\tau(\bar{\mu}_{(k+1)}),\bar{\tau}^{-1}(g_{(k)}^{-1}g_{(k+1)}), \nonumber \\
    &Ad^*_{(g_{(k)}^{-1}g_{(k+1)})^{-1}}(\bar{\mu}_{(k+1)})-\bar{\mu}_{(k)}). \label{constraint adapted discretization inverse}
\end{align}
Here $\bar{\mu}$ denotes the component of $\mu$ in $\mathfrak{d}^*$, while $\bar{\tau}^{-1}(g_{(k)}^{-1}g_{(k+1)})$ denotes the component of $\tau^{-1}(g_{(k)}^{-1}g_{(k+1)})$ in $\mathfrak{d}$. 
\subsection{Constraint-Preserving Integrator}
\label{subsec: constriant preserving integrator}
On $G = \mathbb{SO}(3)$, two commonly used choices for the local diffeomorphism $\tau$ are the exponential and the Cayley map, see \cite{Iserles2000}. Applying the seed formula \eqref{integrator on G}, lifted and adapted as in \eqref{integrator on D*}, to the left-trivialized dynamics given by \eqref{Lie Poisson Suslov vector field}, we obtain the following.
\subsubsection{Exponential Map}
Using the exponential map for $\tau$ in \eqref{constraint adapted discretization inverse}, we obtain the following integrator using the expression for the left-logarithmic derivative given in \cite{Iserles2000}:
\begin{equation} \label{exponential integrator}
    \boxed{ \scriptsize
    \begin{aligned}
         &R_{(k+1)} = R_{(k)} \exp{(t\hat{\Omega}_{(k)})} \\
         &\left( I + \frac{1-\cos{(\|t\Omega_{(k)}\|)}}{\|t\Omega_{(k)}\|^2} t\hat{\Omega}_{(k)} + \frac{\|t\Omega_{(k)}\| - \sin{(\|t\Omega_{(k)}\|)}}{\|t\Omega_{(k)}\|^3} t^2 \hat{\Omega}_{(k)}^2 \right) \Pi_{(k)} \\
         &= t \mathcal{I} \Omega_{(k)} \\
        &\bar{\Pi}_{(k+1)} = \exp{(t\hat{\Omega}_{(k)})} \bar{\Pi}_{(k)}.
    \end{aligned}
    }
\end{equation}
Here $\Omega_{(k)}=(\Omega_{1(k)},\Omega_{2(k)},0)$ and $\bar{\Pi}_{(k)} = (\Pi_{1(k)},\Pi_{2(k)},0)$. The matrix $R \in \mathbb{SO}(3)$ represents the orientation of the rigid body, while $\Pi \in \mathbb{R}^3 \cong \mathfrak{so}(3)$ denotes the (body) angular momentum. 
%
%
%
\subsubsection{Cayley Map}
\begin{figure}[h]
    \centering
    \includegraphics[scale=0.34]{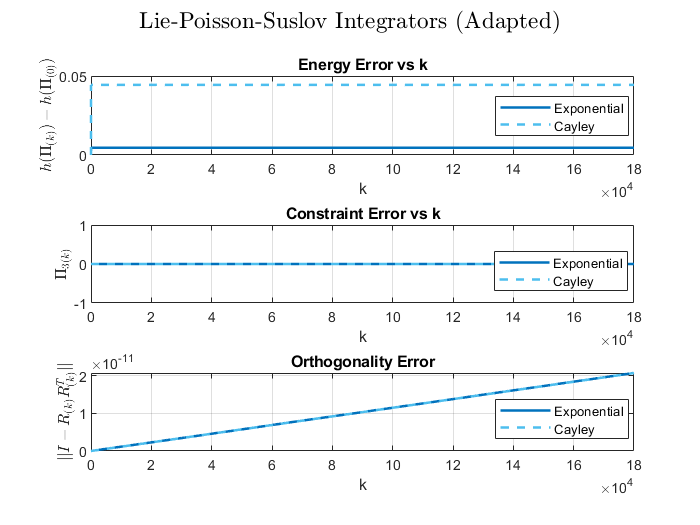}
    \caption{Exponential and Cayley map-based constraint-adapted (Lie-Poisson-Suslov) integrators for the Suslov problem, run for 30 minutes with inertia matrix $\mathcal{I}=diag(1,10,100)$, initial conditions $\Pi_{(0)}=(1,1,0)$, $R_{(0)}=I$, and step size $t =0.01$ s (all in standard units).}
    \label{fig: integrator LPS comparison}
\end{figure}

\begin{figure}[h]
    \centering
    \includegraphics[scale=0.34]{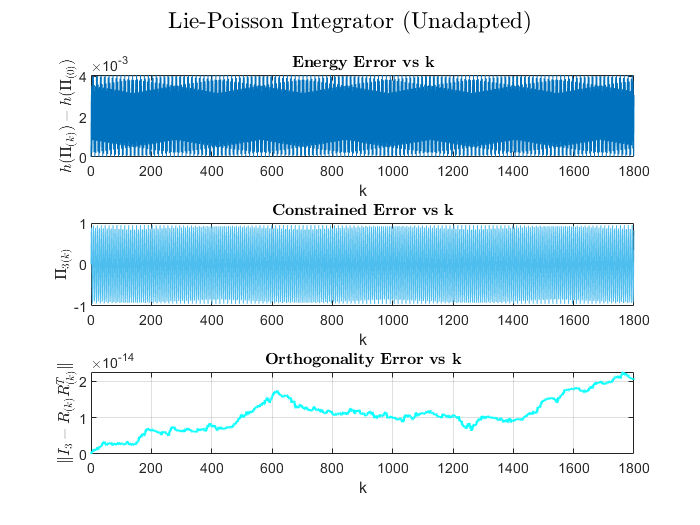}
    \caption{Exponential map-based (Lie-Poisson) integrator for the Suslov problem, run for just 18 seconds with inertia matrix $\mathcal{I}=diag(1,10,100)$, initial conditions $\Pi_{(0)}=(1,1,0)$, $R_{(0)}=I$, and step size $t =0.01$ s (all in standard units).}
    \label{fig: integrator LP}
\end{figure}

Using the Cayley map for $\tau$ in \eqref{constraint adapted discretization inverse}, defined by $Cay(\hat{x}) = (I-\hat{x})^{-1}(I+\hat{x})$, we obtain the following integrator using the expression for the left-logarithmic derivative given in 
\cite{Iserles2000}:
\begin{equation} \label{cayley integrator}
    \boxed{ \scriptsize
    \begin{aligned}
         &R_{(k+1)} = R_{(k)} (I-t\hat{\Omega}_{(k)})^{-1}(I+t\hat{\Omega}_{(k)}) \\
         &2(I+t\hat{\Omega}_{(k)}) \Pi_{(k)} = (1+\|t\Omega_{(k)}\|^2)t \mathcal{I} \Omega_{(k)} \\
        &\bar{\Pi}_{(k+1)} =  (I-t\hat{\Omega}_{(k)})^{-1}(I+t\hat{\Omega}_{(k)})\bar{\Pi}_{(k)}.
    \end{aligned}
    }
\end{equation}
Notice that both integrators \ref{exponential integrator} and \ref{cayley integrator} satisfy the constraint, remain on the Lie group, and exhibit stabilized energy error, as illustrated in Fig. \ref{fig: integrator LPS comparison}. For comparison, the unadapted integrator \eqref{integrator on T*G} is also shown in Fig. \ref{fig: integrator LP}, which clearly fails to preserve the constraint.

\end{document}